\documentclass{article}
\usepackage{amsmath, amsthm, amsfonts, amssymb}
\usepackage{hyperref}
\usepackage[latin1]{inputenc}

\makeatletter \@addtoreset{equation}{section}

\makeatletter

\def\beg   {\begin{theorem}}   \def\ee   {\end{theorem}}
\def\be   {\begin{equation}}   \def\ee   {\end{equation}}
\def\ba   {\begin{array}}      \def\ea   {\end{array}}
\def\bea  {\begin{eqnarray}}   \def\eea  {\end{eqnarray}}
\def\bean {\begin{eqnarray*}}  \def\eean {\end{eqnarray*}}

\newtheorem{lemma}{Lemma}[section]
\newtheorem{theorem} [lemma]{Theorem}

\newtheorem{example}[lemma] {Example}

\newcommand{\mC}{\ensuremath{\mathbb{C}}}
\newcommand{\mD}{\ensuremath{\mathbb{D}}}

\newcommand{\mN}{\ensuremath{\mathbb{N}}}
\newcommand{\mR}{\ensuremath{\mathbb{R}}}

\begin{document}

\vspace{4cm}
\begin{center} \LARGE{\textbf{ A Value Distribution Result and Some Normality Criteria using Partial Sharing of Small Functions }}
\end{center}
\vspace{1cm}

\begin{center} \bf{K. S. Charak$^{1 }$, \quad Shittal Sharma$^{2}$ }
\end{center}

\begin{center}
 Department of Mathematics, University of Jammu,
Jammu-180 006, INDIA.\\
{$^{1}$ E-mail: kscharak7@rediffmail.com }\\

{$^{2}$ E-mail: shittalsharma\_mat07@rediffmail.com }
\end{center}

\bigskip
\begin{abstract}
In this paper, we first generalize a value distribution result of Lahiri and Dewan \cite{IS} and as an application of this result we prove a normality criterion using partial sharing of small functions. Further, in sequel normality criteria of Hu and Meng \cite{HM1} and Ding, Ding and Yuan \cite{DDY} are improved and generalized when the domain $D:=\{z:|z|<R, 0<R\leq\infty \}$.
\end{abstract}

\vspace{1cm}\noindent \textbf{Keywords: } Normal Families, Meromorphic
Functions, Differential Monomials, Sharing of values.

\vspace{0.5cm} \noindent\textbf{AMS subject classification: 30D35, 30D45}

\vspace{4cm}

\normalsize
\newpage

\section{Introduction and Main Results}
We assume that the reader is familiar with the theory of normal
families of meromorphic functions on a domain $D\subseteq \mC ,$
one may refer to \cite{JS}.

\medskip

The idea of sharing of values was introduced in the study of
normality of families of meromorphic functions,
for the first time, by W. Schwick \cite{sch} in 1989.\\
Two non-constant meromorphic functions $f$ and $g$ are said to
share a value $\omega\in \mC$ IM(Ignoring multiplicities) if $f$ and
$g$ have the same $\omega-$points counted with ignoring multiplicities. If
multiplicities of $\omega -$points of $f$ and $g$ are counted, then $f$ and $g$ are said to share the value $\omega$ CM. For deeper insight into
the sharing of values by meromorphic functions, one may refer to \cite{YY}.

\medskip
In this paper all meromorphic functions are considered on $D:=\{z:|z|<R, 0<R\leq \infty \}$ excepting Theorem A and Theorem \ref{THM1}, where the domain is the whole complex plane. A meromorphic function $\omega(z)$ is said to be a {\it small function} of a meromorphic function $f(z)$ if $T(r,\omega)=o\left(T(r,f)\right) \text{ as } r \longrightarrow R.$ Further, we say that a meromorphic function $f$ share a small function $\omega$ {\it
partially} with a meromorphic function $g$ if
$$\overline{E}(\omega,f)=\{z\in \mC: f(z)-\omega(z)=0 \} \subseteq \overline{E}(\omega,g)=\{z\in
\mC: g(z)-\omega(z)=0 \},$$ where $\overline{E}(\omega,\phi)$ denotes the set of
zeros of $\phi - \omega$ counted with ignoring multiplicities.

\medskip

The function of the form $ M[f]= f^{n_0}(f')^{n_{1}} \cdots (f^{(k)})^{n_{k}} $ is called a { \it differential monomial} of $f$ of {\it degree} $d= n_0+n_{1}+\cdots +n_{k}$,
where $n_0, n_{1},\cdots , n_{k}$ are non-negative integers.\\
 
 \medskip
 
In the present discussion, we have used the idea of partial
sharing of small functions in the study of normality of families of
meromorphic functions. One can verify that a good amount of
results on normal families proved by using the sharing of values
can be proved under weaker hypothesis of partial sharing of values or small functions.\\

Lahiri and Dewan \cite{IS} proved the following result:\\
{\bf Theorem A} Let $f$ be a transcendental meromorphic function and $F= (f)^{n_0}(f^{(k)})^{n_1}$, where $n_0 (\geq 2), n_1$ and $k$ are positive integers such that $n_0(n_0-1)+(1+k)(n_0n_1-n_0-n_1)>0.$ Then
$$\left[1-\dfrac{1+k}{n_0+k}-\dfrac{n_0(1+k)}{(n_0+k)\{n_0+(1+k)n_1\}}\right]T(r,F)\leq \overline{N} \left(r,\frac{1}{F-\omega} \right)+S(r,F)$$
for any small function $\omega(\not\equiv 0, \infty)$ of $f$.

\medskip
This is natural to ask whether Theorem A remains valid for a general class of monomials. In this direction, we have proved that it does hold for a larger class of monomials. Precisely, we have
\begin{theorem}
Let $f$ be a transcendental meromorphic function.
Let 
\begin{equation}
F = f^{n_{0}}(f')^{n_{1}} \cdots (f^{(k)})^{n_{k}},
\end{equation} 
where $k,n_{0},n_{1}, \cdots ,n_{k}$ are non-negative
integers with $k\geq 1, n_{0}\geq 2$ and $ n_{k}\geq 1$ such that
\begin{equation}
n_{0}(n_{0}-1) + \sum^{k}_{j=1}(j+1)(n_{0}n_{j}-n_{j}-n_{0})+(k-1)n_{0}>0.
\end{equation}
 Then
 \begin{equation}
\left[1-\dfrac{1+\frac{k(k+1)}{2}}{n_{0}+\frac{k(k+1)}{2}}-\dfrac{n_{0}(1+\frac{k(k+1)}{2})}{\{n_{0}+\frac{k(k+1)}{2}\} \{n_{0}+ \sum^{k}_{j=1}(j+1)n_{j}\}}+o(1)\right]T(r, F)$$
$$ \leq \overline{N} \left(r,\frac{1}{F-\omega} \right)+ S(r, F)
\end{equation}
 for any small function $\omega(\not\equiv 0,\infty)$ of $f.$
 \label{THM1}
\end{theorem}
{\bf Note:} When $f$ has no poles then Theorem \ref{THM1} holds without the condition $(1.2)$.

\medskip
 
 As an application of Theorem \ref{THM1}, we prove a normality criterion using the idea of partial sharing of small functions.
 
 \begin{theorem}
 Let $\mathcal{F}$ be a family of meromorphic functions  such that each $f \in \mathcal{F}$ has only zeros of multiplicity at least $k \geq 2.$ Let $n_{0}, n_{1}, \cdots , n_{k}$ be non- negative integers with $n_{0}\geq 2, n_{k}\geq 1$ such that
$$ n_{0}(n_{0}-1)+\sum ^{k}_{j=1}(j+1)(n_{0}n_{j}-n_{0}-n_{j})+(k-1)n_{0}>0$$
Let $\omega(z)$ be a small function of each $f\in \mathcal{F}$ having no zeros and poles at the origin. If there exists $\widetilde{f} \in \mathcal{F}$ such that $M[f]$ share $\omega$ partially with $M[\widetilde{f}],$ for every $f \in \mathcal{F}$, then $\mathcal{F}$ is a normal family.
 \label{THM2}
 \end{theorem}
Further, one can see that Theorem 4.1 of Hu and Meng \cite{HM1} may be generalized to a class of monomials as
\begin{theorem}
Let $k\in \mN$ and $\mathcal{F}$ be a family of non-constant meromorphic functions such that each $f\in \mathcal{F}$ has only zeros of multiplicity at least $k.$ Let $n_0, n_1, \cdots, n_k$ be non-negative integers with $n_0\geq 2, n_k\geq 1$ such that
$$n_0(n_0-1)+\sum^{k}_{j=1}(j+1)(n_0 n_j -n_0-n_j)+(k-1)n_0>0.$$
Let $\omega(z)$ be a small function of each $f \in \mathcal{F}$ having no zeros and poles at the origin. If, for each $f\in \mathcal{F}, \left(M[f]-\omega \right)(z)=0$ implies $|f^{(k)}(z)|\leq A,$ for some $A>0,$ then $\mathcal{F}$ is a normal family.
\label{THM}
\end{theorem}
 
 \section{Proof of Main Results}

 {\em Proof of Theorem \ref{THM1}:} Since (see \cite{AP}) 
 $$T(r,f)+S(r,f) \leq C T(r,F)+ S(r, F)$$
 and
 $$T(r,F) \leq \left[n_{0}+ \sum^{k}_{j=1}(j+1)n_{j}\right]T(r,f)+S(r,f),$$
 where $C$ is a constant, it follows that $T(r,\omega)=S(r,F)$ as $r\longrightarrow \infty.$ Precisely, $\omega$ is a small function of $f$ iff $\omega$ is a small function of $F$.\\
 Now, by Second Fundamental Theorem of Nevanlinna for three small functions(see \cite{WK} pp. 47), we have
\begin{equation}
[1+o(1)]T(r,F) \leq \overline{N}(r,F)+\overline{N}(r,\frac{1}{F})+ \overline{N} \left(r,\frac{1}{F-\omega} \right)+ S(r,F).
\end{equation}
 Next, we have
 \begin{align*} 
 \overline{N}(r,\frac{1}{F}) & \leq \overline{N}(r,\frac{1}{f})+ \sum^{k}_{j=1} \overline{N_{0}}(r,\frac{1}{f^{(j)}})\\
 & \leq \overline{N}(r,\frac{1}{f})+\sum^{k}_{j=1}j \left[\overline{N}(r,\frac{1}{f})+ \overline{N}(r,f)\right]+S(r,f)\\
 & = \overline{N}(r,\frac{1}{f})+\frac{k(k+1)}{2} \left[\overline{N}(r,\frac{1}{f})+ \overline{N}(r,f)\right]+ S(r,f),
 \end{align*}
 where $ \overline{N} _{0}(r,\frac{1}{f^{(j)}})$ is the number of those zeros of $f^{(j)}$ in $|z|\leq r$ which are not the zeros of $f$.\\
 That is,
\begin{equation}
\overline{N}(r,\frac{1}{F}) \leq \left[1+\frac{k(k+1)}{2}\right] \overline{N}(r,\frac{1}{f})+\frac{k(k+1)}{2}\overline{N}(r,f)+ S(r,f).
 \end{equation}
 Also, we can see that 
\begin{equation}
N(r,\frac{1}{F})-\overline{N}(r,\frac{1}{F}) \geq \left[(k+1)n_{0}+ \sum ^{k}_{j=1}n_{j}-1\right] \overline{N}_{(k+1}(r,\frac{1}{f})+(n_{0}-1) \overline{N}_{k)}(r,\frac{1}{f}),
\end{equation}
 
where $\overline{N}_{(k+1}(r,\frac{1}{f})$ and $\overline{N}_{k)}(r,\frac{1}{f})$ are the counting functions ignoring multiplicities of those zeros of $f$ whose multiplicity is $\geq k+1$ and $\leq k$ respectively.\\
 
 Now from (2.2) and (2.3), we get
 \begin{align*}
 \overline{N}(r,\frac{1}{F}) 
 & \leq \left[1+\frac{k(k+1)}{2}\right] \overline{N}_{(k+1}(r,\frac{1}{f})\\
 & +\dfrac{\left[1+\frac{k(k+1)}{2}\right]}{n_{0}-1} \left[ N(r,\frac{1}{F}) -\overline{N}(r,\frac{1}{F})- \left((k+1)n_{0}+ \sum ^{k}_{j=1}n_{j}-1\right)\overline{N}_{(k+1}(r,\frac{1}{f}) \right]\\
 & + \frac{k(k+1)}{2} \overline{N}(r,f)+ S(r,f).
 \end{align*}
 That is,
 \begin{align*}
 \left[1+\frac{\left(1+\frac{k(k+1)}{2} \right)}{n_{0}-1}\right]\overline{N}(r,\frac{1}{F})
  & \leq \left(1+\frac{k(k+1)}{2}\right) \left(1-\frac{(k+1)n_{0}+ \sum ^{k}_{j=1}n_{j}-1}{n_{0}-1}\right) \overline{N}_{(k+1}(r,\frac{1}{f})\\
  & +\dfrac{1+\frac{k(k+1)}{2}}{n_{0}-1}N(r,\frac{1}{F})+\dfrac{k(k+1)}{2} \overline{N}(r,f)+ S(r,f).\\
 \end{align*}
 
 Since $\overline{N}(r,f)=\overline{N}(r,F)$ and $S(r,f)= S(r,F),$ we have 
 \begin{align*}
\overline{N}(r,\frac{1}{F}) & \leq \dfrac{1+\frac{k(k+1)}{2}}{n_{0}+\frac{k(k+1)}{2}}N(r,\frac{1}{F})+ \dfrac{(\frac{k(k+1)}{2})(n_{0}-1)}{n_{0}+\frac{k(k+1)}{2}}\overline{N}(r,f)+ S(r,f)\\
 & = \dfrac{1+\frac{k(k+1)}{2}}{n_{0}+\frac{k(k+1)}{2}}N(r,\frac{1}{F})+ \dfrac{(\frac{k(k+1)}{2})(n_{0}-1)}{n_{0}+\frac{k(k+1)}{2}}\overline{N}(r,F)+ S(r,F).
 \end{align*}
 
 Therefore, (2.1) yields  
 
 \begin{equation}
 [1+o(1)]T(r,F) \leq \overline{N}\left(r,\frac{1}{F-\omega} \right) + \dfrac{1+\frac{k(k+1)}{2}}{n_{0}+\frac{k(k+1)}{2}}N(r,\frac{1}{F})+\dfrac{n_{0}(1+\frac{k(k+1)}{2})}{n_{0}+\frac{k(k+1)}{2}} \overline{N}(r,F)+S(r,F).
 \end{equation}
 
Also, if $f$ has a pole of multiplicity p, then $F$ has a pole of multiplicity 
$$n_{0}p+n_{1}(p+1)+\cdots +n_{k}(p+k) \geq n_{0}+2n_{1}+\cdots+(k+1)n_{k}=n_{0}+ \sum^{k}_{j=1}(j+1)n_{j}$$
and therefore,

\begin{equation}
N(r,F)\geq \left[n_{0}+ \sum^{k}_{j=1}(j+1)n_{j}\right]\overline{N}(r,F).
\end{equation}

Finally, from (2.4) and (2.5), we find that

\medskip
$[1+o(1)]T(r,F) \leq \overline{N}\left(r,\dfrac{1}{F-\omega} \right) +\dfrac{1+\frac{k(k+1)}{2}}{n_{0}+\frac{k(k+1)}{2}}N \left(r,\dfrac{1}{F} \right)$ $$+\dfrac{n_{0}(1+\frac{k(k+1)}{2})}{(n_{0}+\frac{k(k+1)}{2})(n_{0}+ \sum^{k}_{j=1}(j+1)n_{j})}N(r,F)+S(r,F).$$
That is,\\
$$ \left[1-\dfrac{1+\frac{k(k+1)}{2}}{n_{0}+\frac{k(k+1)}{2}}-\dfrac{n_{0}(1+\frac{k(k+1)}{2})}{(n_{0}+\frac{k(k+1)}{2})(n_{0}+ \sum^{k}_{j=1}(j+1)n_{j})}+o(1) \right]T(r,F)$$ $$ \leq \overline{N}\left(r,\frac{1}{F-\omega} \right) + S(r,F).$$
~~~~~~~~~~~~~~~~~~~~~~~~~~~~~~~~~~~~~~~~~~~~~~~~~~~~~~~~~~~~~~~~~~~~~~~~~~~~~~~~~~~~~~~~~~~~~~~~~~~~~~~~~~~~~~~~~~~~~~~~~~~~~~~~~~~~~~~$\Box$
For the proof of Theorem \ref{THM2}, besides Theorem \ref{THM1}, we also need the following lemma which is a straight forward generalization of $Lemma$ 3 in \cite{DDY}.

\begin{lemma}
Let $f$ be a non-constant rational function with only zeros of multiplicity at least $k$, where $k \geq 2.$ Let $n_0, n_1, n_2, \cdots, n_k$ be non-negative integers with $n_0 \geq 2$ and $n_k \geq 1.$ Let $\omega \neq 0$ be a finite complex number. Then $M[f]-\omega$ has at least two distinct zeros.\label{lemma5}
\end{lemma}

{\em Proof of Theorem \ref{THM2}:} Since normality is a local property, we may assume that $D=\mD.$ Suppose $\mathcal{F}$ is not normal in $\mD$. 
In particular, suppose that $\mathcal{F}$ is not normal at $z=0.$ Then, by Zalcman's lemma (see \cite{LZ}), there exist  a sequence $\{f_{n}\}$ of functions in $\mathcal{F}$, a sequence $\{ z_{n} \}$ of complex numbers in $\mD$ with $z_{n}\longrightarrow 0$ as $n\longrightarrow \infty,$ and  a sequence $ \{ \rho_{n} \} $ of positive real numbers with $\rho_n \longrightarrow 0$ as $n \longrightarrow \infty$ such that the sequence $\{g_{n}\}$ defined by
 $$g_{n}(z)=\rho^{-\alpha}f_{n}(z_{n}+\rho_{n}z); 0\leq \alpha < k,$$
 converges locally uniformly to a non-constant meromorphic function $g(z)$ in $\mC$ with respect to the spherical metric. Moreover, $g(z)$ is of order at most 2. By Hurwitz's theorem, the zeros of $g(z)$ have multiplicity at least $k.$\\
 Let $\alpha = \dfrac{\sum^{k}_{j=1}jn_{j}}{\sum^{k}_{j=0}n_{j}}<k.$ Then 
\begin{align*} 
M[g_n](z) & = \left(g_n(z)\right)^{n_o}\left(g'_n(z)\right)^{n_1}\cdots \left(g^{(k)}_n(z)\right)^{n_k}\\
& = \rho_{n}^{-\alpha n_{0}}\left(f_n(z_n+\rho_n z)\right)^{n_0} \rho_{n}^{-\alpha n_1+n_1}\left(f'_n(z_n+\rho_n z)\right)^{n_1} \cdots \rho_{n}^{-\alpha n_k+kn_k}\left(f^{(k)}_n(z_n+\rho_n z)\right)^{n_k}\\
& = \rho^{-\alpha \sum^{k}_{j=0}n_j + \sum^{k}_{j=1}jn_j}\left(f_n(z_n+\rho_n z)\right)^{n_0}\left(f'_n(z_n+\rho_n z)\right)^{n_1}\cdots\left(f^{(k)}_n(z_n+\rho_n z)\right)^{n_k}\\
& =M[f_n](z_n+\rho_n z).
\end{align*}
 On every compact subset of $\mC$ that contains no poles of $g$, we have
 $$M[f_n](z_n+\rho_nz)-\omega(z_n+\rho_nz)=M[g_n](z)-\omega(z_n+\rho_nz)\longrightarrow M[g](z)-\omega_0$$
 spherically uniformly, where $\omega_0= \omega(0).$\\
 Since $g$ is a non-constant meromorphic function of order at most 2 and $\omega_0\neq 0, \infty,$ it immediately follows that $M[g]\not\equiv \omega_0.$ Using Theorem \ref{THM1} and Lemma \ref{lemma5}, $M[g]-\omega_0$ has at least two distinct zeros, say, $w_0$ and $v_0$. Choose $r>0$ such that the open disks $D(w_0,r)=\{z:|z-w_0|<r\}$ and $ D(v_0,r)=\{z:|z-v_0|<r\}$ are disjoint and their union contains no zeros of $M[g]-\omega_0$ different from $w_0$ and $v_0$ respectively. Then, by Hurwitz's theorem, we see that for sufficiently large $n,$ there exist points $w_{n}\in D(w_0,r)$ and $ v_n \in D(v_0,r)$ such that\\
 $$\left(M[f_n]-\omega\right)(z_n+\rho_nw_n)=0,$$
 and
 $$\left(M[f_n]-\omega\right)(z_n+\rho_nv_n)=0.$$ 
 Since by hypothesis, $M[f_n]$ share $\omega$ partially with $M[\widetilde{f}]$, for every $n,$ it follows that
 $$\left(M[\widetilde{f}]-\omega \right)(z_n+\rho_nw_n)=0,$$
 and
 $$\left(M[\widetilde{f}]-\omega \right)(z_n+\rho_nv_n)=0.$$
 By letting $n\longrightarrow\infty,$ and noting that $z_n+\rho_nw_n\longrightarrow 0, z_n+\rho_nv_n\longrightarrow 0,$ we find that
 $$\left(M[\widetilde{f}]-\omega \right)(0)=0.$$
 Since the zeros of $M[\widetilde{f}]-\omega$ have no accumulation point, $z_n+\rho_nw_n=0$ and $ z_n+\rho_nv_n=0$ for sufficiently large $n$. That is, $D(w_0,r)\cap D(v_0,r)\neq \phi ,$ a contradiction.~~~~~~~~~~~~~~~~~~~~~~~~~~~~~~~~~~~~~~~~~~~~~~~~~~~~~~~~~~~~~~~~~~~~~~~~~~~~~~~~~~~~~~~~~~~~~~~~~$\Box$

\bigskip

{\em Proof of Theorem \ref{THM}:} As established in the proof of Theorem \ref{THM2}, we similarly find that $M[g]\not\equiv \omega_0.$ By Theorem \ref{THM1} and Lemma 2.6 in \cite{ZXY}, $M[g]-\omega_0$ has at least one zero $w_0,$ say. By Hurwitz's Theorem, there is a sequence of complex numbers $\{w_n \}$ such that $w_n \longrightarrow w_0$ as $n\longrightarrow \infty,$ and
$$\left(M[f_n]-\omega \right)(z_n+\rho_n w_n)=0$$
Again, since $k>\alpha,$
\begin{align*}
|g^{(k)}_{n}(w_n)| & = \rho^{k-\alpha}_{n}|f^{(k)}_{n}(z_n+\rho_n w_n)|\\
 & \leq \rho^{(k-\alpha)}_{n} A\\
 & = A \rho^{k-\frac{\sum^{k}_{j=1}jn_j}{\sum^{k}_{j-0}n_j}}_n \longrightarrow 0 \text{ as } n \longrightarrow \infty.
\end{align*}
Therefore, $g^{(k)}(w_0)= \lim_{n\longrightarrow \infty}g^{(k)}_n(w_n)=0$\\
$\Rightarrow M[g](w_0)=0 \neq \omega_0,$ which is a contradiction.~~~~~~~~~~~~~~~~~~~~~~~~~~~~~~~~~~~~~~~~~~~~~~~~~~~~~~~~~~~~~~~$\Box$

\section{ Conclusions}
  Though our results do generalize and improve the results of Hu and Meng \cite{HM1} and Ding, Ding and Yuan\cite{DDY} when the domain $D$ is $\{z:|z|<R, \ 0,R\leq \infty \}$, there seems no way of proving our results on arbitrary domain since the idea of small function on arbitrary domain is not available, as for as we know. However, by making certain modifications in the proofs of results of Hu and Meng\cite{HM1} and Ding, Ding and Yuan\cite{DDY}, one can easily extend and improve these results on arbitrary domain with shared value being a non-zero complex value. Precisely, one obtains,
  
\begin{theorem} 
Let $\mathcal{F}$ be a family of non-constant meromorphic
functions on a domain $D$  with all zeros of each $f \in
\mathcal{F}$ having multiplicity at least $k$, where $k \geq 2$.
Let $\omega \neq 0$ be a finite complex number and $n_0, n_{1},\cdots
,n_{k}$ be non-negative integers with $n_0 \geq 2$ and $n_{1}+n_{2}
\cdots +n_{k} \geq 1.$ If there exists $\widetilde{f} \in \mathcal{F}
$ such that $M[f]$ share $\omega$ partially with $M[\widetilde{f}]$ for
every $f \in \mathcal{F}$, then $\mathcal{F}$ is  normal on $D$.
\label{THM3}
\end{theorem}

The condition that $f$ has only zeros of multiplicity atleast $k$
in Theorem \ref{THM3} is sharp. For example, consider the open unit disk $\mathbb{D}$, an integer $k \geq 2$, a non-zero complex number $\omega$ and the family
$$\mathcal{F} = \{f_{m}(z)=mz^{k-1}; m=1,2,3,\cdots\}$$
Obviously, each $f_{m} \in \mathcal{F}$ has only a zero of multiplicity $k-1$, and for distinct positive integers $m$, and $l$; we find that $f_{m}^{2}f_{m}^{(k)}$ and $f_{l}^{2}f_{l}^{(k)}$ share $\omega$ IM and $\mathcal{F}$ is not normal at $z=0.$

\medskip

Also, $\omega \neq 0$ in Theorem \ref{THM3} is essential. For example, let $\mathcal{F} = \{f_{^m}\}$, where $f_{m}(z) =
\frac{1}{e^{mz}+1}$; $ m = 1, 2,\cdots $ and $z \in \mathbb{D}$.
Choose $k = 2$, $n = 2$, $n_{1} = 1$, and $n_{2} = 0$, we have
$$M[f_{m}] = f_{m}^{2}f'_{m} = -\frac{me^{mz}}{(e^{mz}+1)^{4}} \neq 0.$$
Thus, for any $f, g \in \mathcal{F}$, $M[f]$ and $M[g]$ share $0$
IM. But we see that $\mathcal{F}$ is not normal in $\mathbb{D}$.

\medskip

\begin{theorem}
Let $\mathcal{F}$ be a family of non-constant holomorphic
functions on a domain $D$ with all zeros of each $f \in
\mathcal{F}$ having multiplicity at least $k$, where $k \geq 2$.
Let $\omega \neq 0$ be a finite complex number and $n_0, n_{1},\cdots
,n_{k}$ be non-negative integers with $n_0 \geq 1$ and $n_{1}+n_{2}
\cdots +n_{k} \geq 1.$ If there exists $\widetilde{f} \in \mathcal{F}
$ such that $M[f]$ share $\omega$ partially with $M[\widetilde{f}]$ for
every $f \in \mathcal{F}$, then $\mathcal{F}$ is normal on $D$.
\label{THM4}
\end{theorem}

As an illustration of Theorem \ref{THM4}, we have the following example:
\begin{example}
Consider $\mathcal{F}=\{f_{m}(z)=me^{\frac{z}{m}}: m\in \mathbb{N} \},$
defined on $\mC.$ Take $k=2, n=1, n_{1}=0,$ and $n_{2}=1.$ Then
$$M[f_{m}]= f_{m}f''_{m}= e^{\frac{2z}{m}},$$
and $M[f_{m}]=1$ iff $\dfrac{2z}{m}=2k\pi i, k\in \mathbb{Z}$ iff $z=mk\pi i$\\
For
$$m=1; z=0, \pm \pi i,  \pm2\pi i, \pm3\pi i, \cdots$$
$$m=2; z=0, \pm2\pi i,  \pm4\pi i, \pm6\pi i, \cdots$$
$$m=3; z=0, \pm3\pi i,  \pm6\pi i, \pm9\pi i, \cdots$$
and so on.\\
Thus for each $m \geq 2$, $M[f_{m}]$ share 1 partially with
$M[f_{1}]$. Next, we have $\forall z, |z|\leq r, r>0;
|f_{m}(z)|=|me^{\frac{z}{m}}|=me^{\frac{\Re(z)}{m}}<
me^{\frac{r}{m}}=M, $ say, where $M>0$ depends on $r$ and this is
true for each $m \in \mathbb{N}$. That is, $\mathcal{F}$ is locally
bounded on $\mC$ and hence by Montel Theorem $\mathcal{F}$ is
normal.
\end{example}

\bibliographystyle{amsplain}

\end{document}